\documentclass[amstex,12pt,russian,amssymb]{article}

\usepackage{mathtext}
\usepackage[cp1251]{inputenc}
\usepackage[T2A]{fontenc}
\usepackage[russian]{babel}
\usepackage[dvips]{graphicx}
\usepackage{amsmath}
\usepackage{amssymb}
\usepackage{amsxtra}
\usepackage{latexsym}
\usepackage{ifthen}

\def\XXint#1#2#3{{\setbox0=\hbox{$#1{#2#3}{\int}$}
\vcenter{\hbox{$#2#3$}}\kern-.5\wd0}}

\begin{document}

\newcounter{lemma}
\newcommand{\lemma}{\par \refstepcounter{lemma}%
{\bf Лемма \arabic{lemma}.}}

\newcounter{corollary}
\newcommand{\corollary}{\par \refstepcounter{corollary}%
{\bf Следствие \arabic{corollary}.}}

\newcounter{remark}
\newcommand{\remark}{\par \refstepcounter{remark}%
{\bf Замечание \arabic{remark}.}}

\newcounter{theorem}
\newcommand{\theorem}{\par \refstepcounter{theorem}%
{\bf Теорема \arabic{theorem}.}}

\newcounter{proposition}
\newcommand{\proposition}{\par \refstepcounter{proposition}%
{\bf Предложение \arabic{proposition}.}}

\renewcommand{\refname}{\centerline{\bf Список литературы}}

\newcommand{\proof}{{\it Доказательство.\,\,}}

\noindent УДК 517.5

{\bf Р.Р.~Салимов} (Институт математики НАН Украины),

{\bf Е.А.~Севостьянов} (Житомирский государственный университет им.\
И.~Франко)

\medskip
{\bf Р.Р.~Салімов} (Інститут математики НАН України),

{\bf Є.О.~Севостьянов} (Житомирський державний університет ім.\
І.~Франко)

\medskip
{\bf R.R.~Salimov} (Institute of Mathematics of NAS of Ukraine),

{\bf E.A.~Sevost'yanov} (Zhitomir State University of I.~Franko)

\medskip
{\bf Об абсолютной непрерывности отображений, искажающих модули
цилиндров}

\medskip
{\bf Про абсолютну неперернвість відображень, що спотворюють модулі
циліндрів}

\medskip
{\bf On absolute continuity for mappings, distorting moduli of
cylinders}

\medskip
В работе рассматриваются отображения, удовлетворяющие одному
модульному неравенству относительно цилиндров в пространстве.
Искажение модуля мажорируется интегралом, зависящим от некоторой
локально интегрируемой функции. Доказан результат об абсолютной
непрерывности на линиях указанных отображений

\medskip
В роботі розглянуто відображення, що задовольняють одну модульну
нерівність відносно циліндрів у просторі. Спотворення модуля
мажорується інтегралом, що залежить від деякої локально інтегровної
функції. Доведено результат про абсолютну неперервність на лініях
вказаних відображень

\medskip
In the present paper, mappings satisfying one modular inequality
with respect to cylinders in a space, are considered. Distorting of
modulus is majorized by an integral which depends from some locally
integrable function. The result on absolute continuity on lines of
the mappings mentioned above is proved

\newpage
{\bf 1. Введение.} Настоящая заметка посвящена изучению свойств
отображений с конечным искажением, активно изучаемых в последнее
время (см., напр., \cite{IM}, \cite{MRSY}--\cite{MRSY$_1$},
\cite{GRSY}, \cite{GG} и \cite{GS}). Здесь и далее $D$ -- область в
${\Bbb R}^n,$ $n\geqslant 2,$ $m$ -- мера Лебега в ${\Bbb R}^n,$
отображение $f:D\rightarrow {\Bbb R}^n$ есть непрерывное
соответствие $f(x)=(f_1(x),\ldots, f_n(x)),$ где $x=(x_1,\ldots,
x_n).$ Напомним, что {\it кривой} $\gamma$ мы называем непрерывное
отображение отрезка $[a,b]$ (открытого, либо полуоткрытого интервала
одного из видов: $(a,b),$ $[a,b),$ $(a,b]$) в ${\Bbb R}^n,$
$\gamma:[a,b]\rightarrow {\Bbb R}^n.$ Кривая будет называться {\it
дугой}, если указанное отображение отрезка (интервала) в ${\Bbb
R}^n$ гомеоморфно (см. \cite[определение~2.7]{Va}). Под семейством
кривых $\Gamma$ подразумевается некоторый фиксированный набор кривых
$\gamma,$ а $f(\Gamma)=\left\{f\circ\gamma|\gamma\in\Gamma\right\}.$
Следующие определения могут быть найдены, напр., в \cite[разд.~1--6,
гл.~I]{Va}.

Борелева функция $\rho:{\Bbb R}^n\,\rightarrow [0,\infty]$
называется {\it допустимой} для семейства $\Gamma$ кривых $\gamma$ в
${\Bbb R}^n,$ если
$\int\limits_{\gamma}\rho (x)\ \ |dx|\geqslant 1$
%
%
для всех кривых $ \gamma \in \Gamma.$ В этом случае мы пишем: $\rho
\in {\rm adm} \,\Gamma .$
{\it Модулем} семейства кривых $\Gamma$ порядка $p\geqslant 1$
называется величина
\begin{equation}\label{eq2**} M_p(\Gamma)=\inf_{\rho \in \,{\rm adm}\,\Gamma} \int\limits_D
\rho^p (x) dm(x)\,.
\end{equation}
Полагаем $M(\Gamma):=M_n(\Gamma).$ Под {\it цилиндром в ${\Bbb
R}^n,$} $n\geqslant 2,$ мы будем понимать следующую тройку: $Z=(Q,
E_1, E_2),$ где $Q$ -- область в ${\Bbb R}^n,$ а $E_1,$ $E_2$ --
подмножества границы $Q$ со следующим свойством: найдётся
гомеоморфизм множества $\overline{Q}$ на единичный цилиндр
$x_1^2+\cdots +x^2_{n-1}\leqslant 1,$ $0\leqslant x_n\leqslant 1,$
отображающий $E_1$ и $E_2$ на его основания (см. \cite{Va$_1$}). С
каждым цилиндром $Z$ будем ассоциировать семейство $\Gamma_Z,$
состоящее из дуг, соединяющих $E_1$ и $E_2$ в $Q.$ Модуль $\Gamma_Z$
будем называть модулем цилиндра $Z$ и обозначать $M(Z).$

\medskip
Пусть $G$ -- открытое множество в ${\Bbb R}^n$ и $I=\{x\in{\Bbb
R}^n:a_i<x_i<b_i,i=1,\ldots,n\}$ -- открытый $n$-мерный интервал.
Отображение $f:I\rightarrow{\Bbb R}^n$ {\it принадлежит классу
$ACL$} ({\it абсолютно непрерывно на линиях}), если $f$ абсолютно
непрерывно на почти всех линейных сегментах в $I,$ параллельных
координатным осям. Отображение $f:G\rightarrow{\Bbb R}^n$ {\it
принадлежит классу $ACL$} в $G,$ когда сужение $f|_I$ принадлежит
классу $ACL$ для каждого интервала $I,$ $\overline{I}\subset G.$

\medskip
В работе \cite{Va$_1$} установлен следующий факт (см. лемму 2).
Пусть $f:D\rightarrow {\Bbb R}^n$ -- гомеоморфизм, удовлетворяющий
условию $M(Z)\leqslant K\cdot M(f(Z))$ для каждого цилиндра
$Z\subset D$ и некоторой постоянной $K\geqslant 1,$ тогда $f\in ACL$
в $D.$ По этому поводу следует также упомянуть интересный результат
Тенгвалля в этом направлении, см. \cite{T}, а также работы авторов,
где установлено свойство $ACL$ для $Q$-отображений, см. \cite{SS}.
Следующее усиление этого результата будет установлено в настоящей
работе.

\medskip
\begin{theorem}\label{th1A}
{\sl Пусть $p>1,$ $Q:{\Bbb R}^n\rightarrow [0, \infty]$ -- локально
интегрируемая функция, $f:D\rightarrow {\Bbb R}^n$ -- открытое
дискретное отображение, удовлетворяющее условию
\begin{equation}\label{eq1}
M_p(Z)\leqslant\int\limits_{f(D)}Q(y)\cdot \rho_*^p(y)dm(y)
\end{equation}
для любого цилиндра $Z\subset D$ и произвольной функции  $\rho_*\in
{\rm adm\,}\Gamma_Z.$ Тогда $f\in ACL.$ }
\end{theorem}

\medskip
{\bf 2. Некоторые предварительные сведения.} Пусть $U$ -- открытое
множество в ${\Bbb R}^n.$ Обозначим через ${\rm Bor\,}U$ класс всех
борелевских подмножеств $U.$ Функция $\varphi:{\rm
Bor\,}U\rightarrow {\Bbb R}$ называется {\it
$q$-ква\-зи\-а\-д\-ди\-тив\-ной}, $q\geqslant 1,$ если выполняются
следующе условия:

1) $\varphi(A)\geqslant 0$ для произвольного борелевского множества
$A\subset U;$

2) из условия $A\subset B$ вытекает неравенство $\varphi(A)\leqslant
\varphi (B)$ какие бы ни были борелевские множества $A, B\subset
{\rm Bor\,} U;$

3) $\varphi(A)<\infty$ для произвольного компактного множества
$A\subset U;$

4) если множества $A_1,\ldots,A_m\subset U$ не пересекаются и
$A_i\subset A\subset U,$ $i=1,\ldots,m,$ то
\begin{equation}\label{eq1.6A}
\sum\limits_{i=1}^m\varphi(A_i)\leqslant q\cdot\varphi(A)\,.
\end{equation}

\medskip
Верхняя и нижняя производные $q$-ква\-зи\-ад\-д\-и\-тив\-ной функции
$\varphi$ в точке $x\in U$ определяются следующим образом:
$$\overline{\varphi^{\,\prime}}(x)=\lim\limits_{h\rightarrow 0}
\sup\limits_{d(Q)<h}\frac{\varphi(Q)}{m(Q)}\,,\quad
\underline{\varphi^{\,\prime}}(x)=\lim\limits_{h\rightarrow 0}
\inf\limits_{d(Q)<h}\frac{\varphi(Q)}{m(Q)}\,,$$
где $Q$ пробегает все открытые кубы и шары такие, что $x\in Q\subset
U.$ Следующее утверждение см., напр., в \cite[лемма~2.3]{MRV$_1$}.

\medskip
\begin{proposition}\label{pr1.1}
{\sl\, Предположим, что функция $\varphi:{\rm Bor\,}U\rightarrow
{\Bbb R}$ является $q$-ква\-зи\-а\-д\-ди\-тив\-ной функцией
множеств. Тогда: 1) функции $\overline{\varphi^{\,\prime}}$ и
$\underline{\varphi^{\,\prime}}$ являются борелевскими, 2) для почти
всех $x\in U$ имеет место неравенство
$\overline{\varphi^{\,\prime}}(x)\leqslant q\cdot
\underline{\varphi^{\,\prime}}(x)<\infty,$
3) для каждого открытого множества $V\subset U,$
$\int\limits_V \underline{\varphi^{\,\prime}}(x)dm(x)\leqslant
q\cdot \varphi(V).$ }
\end{proposition}

Для отображения $f:D\,\rightarrow\,{\Bbb R}^n,$ множества $E\subset
D$ и $y\,\in\,{\Bbb R}^n,$  определим {\it функцию кратности
$N(y,f,E)$} как число прообразов точки $y$ во множестве $E,$ т.е.
$$
N(y,f,E)\,=\,{\rm card}\,\left\{x\in E: f(x)=y\right\}\,,
$$
$$
N(f,E)\,=\,\sup\limits_{y\in{\Bbb R}^n}\,N(y,f,E)\,.
$$
Заметим, что для открытых дискретных отображений
$f:D\,\rightarrow\,{\Bbb R}^n$ всегда $N(f, G)<\infty$ для любой
области $G$ такой, что $\overline{G}\subset D$ (см.
\cite[лемма~2.12(3)]{MRV$_1$}).

\medskip
{\bf 3. Доказательство основного результата.}  Пусть $I$ --
$n$-мерный интервал в ${\Bbb R}^n$ с ребрами параллельными осям
координат и $\overline{I}\subset D.$ Тогда $I=I_0\times J,$ где
$I_0$ -- $(n-1)$-мерный интервал в ${\Bbb R}^{n-1},$ $J$ --
одномерный интервал, $J=(a,b).$ Далее отождествляем ${\Bbb
R}^{n-1}\times {\Bbb R}$ с ${\Bbb R}^n.$ Покажем, что для почти всех
сегментов $J_{z}=\{z\} \times J\,, z\in I_0, $ отображение
$f|_{J_z}$ абсолютно непрерывно. Рассмотрим функцию множеств,
определенную над алгеброй борелевских множеств $B$ в $I_0,$
$$\Phi(B)=\int\limits_{f\left(B\times J\right)}Q(y)dm(y)\,.$$
Заметим, что $\Phi$ -- $q$-квазиаддитивная функция при $q=N(f, I).$
Действительно,
$$\sum\limits_{i=1}^k\Phi(B_i)= \sum\limits_{i=1}^k
\int\limits_{f(B_i\times J)} Q(y)dm(y)\leqslant
$$$$\leqslant\sum\limits_{i=1}^k \int\limits_{f(B\times J)} N(y, f, B_i\times
J)Q(y)dm(y)\leqslant N(f, I)\Phi(B)\,,$$
где $B_i\subset B\subset I_0$ -- борелевские множества, $B_i\cap
B_l=\varnothing$ при $l\ne i.$ По предложению \ref{pr1.1}
\begin{equation}\label{eq11*!}
\varphi(z)=\limsup\limits_{r\rightarrow
0}\frac{\Phi(B(z,r))}{\Omega_{n-1}r^{n-1}}<\infty
\end{equation}
для почти всех $z\in I_0,$ где через $B(z,r)$ обозначается шар в
$I_0\subset{\Bbb R}^{n-1}$ с центром в точке $z\in I_0$ радиуса $r,$
$\Omega_{n-1}$ -- объем единичного  шара в ${\Bbb R}^{n-1}.$
Докажем, что отображение $f$ абсолютно непрерывно на каждом сегменте
$J_z, z\in I_0,$ где предел (\ref{eq11*!}) существует и конечен.
Обозначим соответствующее множество $z$ через $Z_0$ и покажем, что
сумма диаметров образов любого конечного набора непересекающихся
сегментов в $J_z=\{z\}\times J,$ $z\in Z_0,$ стремится к нулю вместе
с суммарной длиной интервалов. Ввиду непрерывности $f$ вдоль $J_z,$
достаточно проверить этот факт для сегментов с рациональными концами
в  $J_z.$ Обозначим $\Delta_i:=(\alpha_i,\beta_i).$

\medskip
Не ограничивая общности, можно считать, что $|f(z, \alpha_i)-f(z,
\beta_i)|\ne 0$ для каждого $i=1,\ldots,k.$ Тогда, поскольку
отображение $f$ непрерывно, для каждого $i=1,\ldots,k $ найдётся
$\delta_i>0$ такое, что
\begin{equation}\label{eq3}
|f(z, \alpha_i)-f(x)|<\frac{|f(z, \alpha_i)-f(z,
\beta_i)|}{4}\,,\quad |x-(z, \alpha_i)|<\delta_i\,,
\end{equation}
\begin{equation}\label{eq4}
|f(z, \beta_i)-f(x)|<\frac{|f(z, \alpha_i)-f(z,
\beta_i)|}{4}\,,\quad |x-(z, \beta_i)|<\delta_i\,.
\end{equation}
В дальнейшем положим $\delta:=\min\limits_{i=1,\ldots, k}\delta_i,$
$0<r<\min\{\delta, {\rm dist}\,(I,
\partial D)\}.$
Множество $\{z\}\times\{\Delta_i\}_{i=1}^k$ покроем цилиндрами
$C_i=B(z, r)\times (\alpha_i, \beta_i)$ и определим
$\widetilde{\rho}_i(z)=2 |f(z, \alpha_i) - f(z, \beta_i)|^{-1}
\cdot\chi_{f(C_i)}(z),$ $i=1,\ldots, k.$ Обозначим также $G_i:=|f(z,
\alpha_i) - f(z, \beta_i)|.$
Определим семейство кривых $\Gamma_i$ следующим образом:
$$\Gamma_i=\{\gamma_x^i: x\in B(z, r/2)\}\,,$$
где $\gamma_x^i:[\alpha_i, \beta_i]\rightarrow {\Bbb R}^n$ --
кривая, определённая как $\gamma_x^i(t)=(x, t).$ Тогда для каждой
локально спрямляемой кривой $f\circ \gamma,$ $\gamma\in \Gamma_i,$
ввиду соотношений (\ref{eq3}) и (\ref{eq4}) будем иметь
$$\int\limits_{f\circ \gamma} \widetilde{\rho}_i(z) |dz|=
2G^{\,-1}_i \int\limits_{f\circ \gamma}\chi_{f(C_i)}(z)|dz|\geqslant
G^{\,-1}_i\cdot |f(z, \alpha_i) - f(z, \beta_i)| = 1\,.$$ В таком
случае,
$\widetilde{\rho}_i\in {\rm adm}\,f(\Gamma_i)$ и, значит, ввиду
условия (\ref{eq1}), получаем:
\begin{equation}\label{eq2}
M_p(\Gamma_i)\leqslant 2^p G_i^{-p}\int\limits_{f(C_i)}Q(y)dm(y)\,.
\end{equation}
Заметим, что
$M_p(\Gamma_i)=\frac{\Omega_{n-1}r^{n-1}}{2^{n-1}|\alpha_i-\beta_i|^{p-1}}$
(см. \cite[пункт~7.2]{Va}), поэтому из  (\ref{eq2}) немедленно
получаем, что
$$\frac{\Omega_{n-1}r^{n-1}}{2^{n-1}|\alpha_i-\beta_i|^{p-1}}\leqslant
2^p G_i^{-p}\int\limits_{f(C_i)}Q(y)dm(y)\,,$$
откуда
\begin{equation}\label{eq5}
|f(z, \alpha_i)-f(z,
\beta_i)|^p\leqslant|\alpha_i-\beta_i|^{p-1}\cdot
\frac{2^{p+n-1}}{\Omega_{n-1}}\cdot\frac{1}{r^{n-1}}\int\limits_{f(C_i)}Q(y)dm(y)\,.
\end{equation}
Обозначим для удобства $V_i:=\int\limits_{f(C_i)}Q(y)dm(y),$ тогда
из (\ref{eq5}) вытекает, что
\begin{equation}\label{eq6}
|\alpha_i-\beta_i|\geqslant \left(\frac{C r^{n-1}|f(z,
\alpha_i)-f(z, \beta_i)|^p}{V_i}\right)^{\frac{1}{p-1}}\,,
\end{equation}
где $C$ -- некоторая постоянная, зависящая только от размерности
пространства $n$ и числа $p.$
Ввиду неравенства Гёльдера
$\frac{\left(\sum\limits_{i=1}^k G_i\right)^p}{\sum\limits_{i=1}^k
V_i}\leqslant \left(\sum\limits_{i=1}^k
\frac{G^{p/(p-1)}_i}{V_i^{1/(p-1)}}\right)^{p-1},$ тогда из
(\ref{eq6}) получаем, что
\begin{equation}\label{eq7}
\frac{C\cdot\left(\sum\limits_{i=1}^k
G_i\right)^p}{\sum\limits_{i=1}^k \frac{V_i}{r^{n-1}}}\leqslant
\left(\sum\limits_{i=1}^k |\alpha_i-\beta_i|\right)^{p-1}\,.
\end{equation}
Ввиду того, что $\Phi$ - квазиаддитивная функция множества, мы будем
иметь
$$\sum\limits_{i=1}^k\frac{V_i}{r^{n-1}}\le N(f, I)\frac{\Phi(B(z, r))}{r^{n-1}}\,.$$
Переходя здесь к верхнему пределу при $r\rightarrow 0,$ учитывая
(\ref{eq11*!}), из (\ref{eq7}) будем иметь
\begin{equation}\label{eq8}
\left(\sum\limits_{i=1}^k |f(z, \alpha_i)-f(z,
\beta_i)|\right)^p\leqslant C_1\cdot\left(\sum\limits_{i=1}^k
|\alpha_i-\beta_i|\right)^{p-1}\,,
\end{equation}
где $C_1$ -- некоторая новая постоянная, зависящая только от
размерности пространства $n$ и функции $Q.$ Пусть теперь
$\varepsilon>0$ и $\Delta_i=(\alpha_i,\beta_i),$ $i=1,\ldots, k$ --
система непересекающихся интервалов в $J$ такая, что
$\sum\limits_{i=1}^k |\alpha_i-\beta_i|<\varepsilon,$ тогда из
(\ref{eq8}) следует, что
$$\left(\sum\limits_{i=1}^k |f(z, \alpha_i)-f(z,
\beta_i)|\right)^p\leqslant C_1\cdot{\varepsilon}^{p-1}\,.$$
Теорема доказана. $\Box$

\medskip
\noindent{{\bf Руслан Радикович Салимов} \\
Институт математики НАН Украины \\
ул. Терещенковская, д. 3 \\
г. Киев-4, Украина, 01 601\\
тел. +38 095 630 85 92 (моб.), e-mail: ruslan623@yandex.ru}

\medskip
\noindent{{\bf Евгений Александрович Севостьянов} \\
Житомирский государственный университет им.\ И.~Франко\\
ул. Большая Бердичевская, 40 \\
г.~Житомир, Украина, 10 008 \\ тел. +38 066 959 50 34 (моб.),
e-mail: esevostyanov2009@mail.ru}

\end{document}